# ASYMPTOTICS FOR SPHERICAL NEEDLETS

By P. Baldi,[1] G. Kerkyacharian,[2] D. Marinucci[1] and D. Picard

*Università di Roma Tor Vergata, LPMA Université de Paris X, Università di Roma Tor Vergata and LPMA Université de Paris 7*

We investigate invariant random fields on the sphere using a new type of spherical wavelets, called needlets. These are compactly supported in frequency and enjoy excellent localization properties in real space, with quasi-exponentially decaying tails. We show that, for random fields on the sphere, the needlet coefficients are asymptotically uncorrelated for any fixed angular distance. This property is used to derive CLT and functional CLT convergence results for polynomial functionals of the needlet coefficients: here the asymptotic theory is considered in the high-frequency sense. Our proposals emerge from strong empirical motivations, especially in connection with the analysis of cosmological data sets.

**1. Introduction.** Over the last two decades, wavelets have emerged as one of the most interesting tools of statistical investigation. In this paper we give an application to the statistical analysis of data sets indexed by the unit sphere $\mathbb{S}^2$. This is motivated mostly by the analysis of the Cosmic Microwave Background radiation (hereafter CMB), currently a very active field of research in astrophysics. Every year hundreds of papers appear in physics journals about CMB and the interest on this topic is going to grow in the next few years when the ESA satellite PLANCK will provide a fresh flow of high-resolution data. Examples of spherical data appear also in other areas of the astrophysical sciences [see Angers and Kim (2005)] or outside astrophysics, that is, brain shape modeling and image analysis [see, e.g., Mardia and Patrangenaru (2005), Dryden (2005) and Dette, Melas and Pepelyshev (2005)].

Received June 2006; revised February 2008.
[1]Supported by MIUR PRIN 2006 project *Stochastic Methods in Finance*.
[2]Supported by ACI *Cosmostat 2004*.

CMB data pose a large amount of challenging statistical problems, for instance, estimation of the correlation structure and of the parameters governing this correlation, testing on the law of the field itself (which is predicted to be Gaussian, or very close to Gaussian, by leading physical models of the Big Bang dynamics), detection of outliers in the observed data (which may signal observations of noncosmological origin, i.e., so-called point sources), testing for isotropy and many others [see Genovese et al. (2004a, 2004b), Marinucci (2004) and Marinucci (2006)].

Random fields on the sphere can be investigated using Fourier developments in spherical harmonics. These methods are, however, difficult to adapt when the data are known only on a portion of the spherical surface. This is actually the case of CMB data, as the observation of this field is missing in the equatorial region, due to the direct radiation from the Milky Way.

In this paper we investigate the statistical properties of the so-called needlets. These are a family of spherical wavelets which were introduced by Narcowich, Petrushev and Ward (2006). Needlets enjoy several properties which are not shared by other spherical wavelets. First they enjoy good localization properties in frequency: needlets are compactly supported in the frequency domain with a bounded support which depends explicitly on a user-chosen parameter. On the other hand, needlets enjoy excellent localization properties in real space, with an exponential decay of the tails (see Figure 2 for a typical graph). See Antoine and Vandergheynst (1999) and Antoine et al. (2002) for a different approach to spherical wavelets.

As a major consequence of the localization property both in the frequency and in the space domain, the needlet coefficients are asymptotically uncorrelated as the frequency tends to $\infty$ for any fixed angular distance. This is the first example of such kind of results for any type of spherical wavelets [see Baldi et al. (2008) for a similar result on the torus]. We use this key property to derive a central limit theorem and a functional central limit theorem for general nonlinear statistics of the wavelets coefficients. We discuss how from these results one can derive, for instance, procedures for testing goodness-of-fit on the angular power spectra.

Let us stress again the great advantage of needlets: their ability (due to localization properties) of dealing with data known only on portions of the spherical surface. We remark also that the needlet construction does not rely on any sort of tangent plane approximation which is typically undertaken to implement wavelets on the sphere.

The plan of the paper is as follows. In Section 2 we describe the construction of needlets, following the approach of Narcowich, Petrushev and Ward (2006). In Section 3 we use them to investigate random fields on the sphere and derive the basic correlation inequality. In Section 4 we recall some classical results on the diagram formula, that are needed in Sections 5 and 6 to derive the main convergence results. In Sections 7 and 8 we discuss statistical applications and the effect of missing observations.

**2. Construction of needlets.** This construction is due to Narcowich, Petrushev and Ward (2006). Its aim is essentially to build a very well-localized tight frame constructed using spherical harmonics, as discussed below. It was recently extended to more general Euclidean settings with fruitful statistical applications [see Kerkyacharian et al. (2007)].

Let us denote by $\mathbb{S}^2$, the unit sphere of $\mathbb{R}^3$. There is a unique positive measure on $\mathbb{S}^2$ which is invariant by rotation, with total mass $4\pi$. This measure will be denoted by $dx$. The following decomposition is well known:

$$\mathbb{L}^2 = \bigoplus_{l=0}^{\infty} \mathcal{H}_l, \tag{1}$$

where $\mathbb{L}^2$ denotes the space of square integrable functions on the sphere with respect to $dx$, and $\mathcal{H}_l$ denotes the vector space of the restriction to $\mathbb{S}^2$ of homogeneous polynomials on $\mathbb{R}^3$, of degree $l$, which are harmonic (i.e., $\Delta P = 0$, where $\Delta$ is the Laplacian on $\mathbb{R}^3$). $\mathcal{H}_l$ is called the space of spherical harmonics of degree $l$ [see Stein and Weiss (1971), Chapter 4; Varshalovich, Moskalev and Khersonskiĭ (1988), Chapter 5] and has dimension $2l+1$. The orthogonal projector on $\mathcal{H}_l$ is given by the kernel operator

$$\forall f \in \mathbb{L}^2 \qquad P_{\mathcal{H}_l} f(x) = \int_{\mathbb{S}^2} L_l(\langle x, y\rangle) f(y)\, dy, \tag{2}$$

where $\langle x, y\rangle$ is the standard scalar product of $\mathbb{R}^3$, and $L_l$ is the Legendre polynomial of degree $l$, defined on $[-1, +1]$, verifying

$$\int_{-1}^{1} L_l(t) L_k(t)\, dt = \frac{2k+1}{8\pi^2} \delta_{l,k},$$

where $\delta_{l,k}$ is the Kronecker symbol. Moreover, by definition of the projection operator,

$$L_l(\langle x, y\rangle) = \sum_{m=-l}^{l} Y_{lm}(x) \overline{Y_{lm}(y)},$$

where the spherical harmonics $Y_{lm}, l = 1, 2, 3, \ldots, m = -l, \ldots, l$, form an orthonormal basis of $\mathcal{H}_l$. For an explicit expression of the functions $Y_{lm}$ [see Varshalovich, Moskalev and Khersonskiĭ (1988), Chapter 5]. Let us point out the reproducing property of the projection operators

$$\int_{\mathbb{S}^2} L_l(\langle x, y\rangle) L_k(\langle y, z\rangle)\, dy = \delta_{l,k} L_l(\langle x, z\rangle). \tag{3}$$

The needlet construction is based on two fundamental steps: Littlewood–Paley decomposition and discretization, which are summarized in the two following subsections.

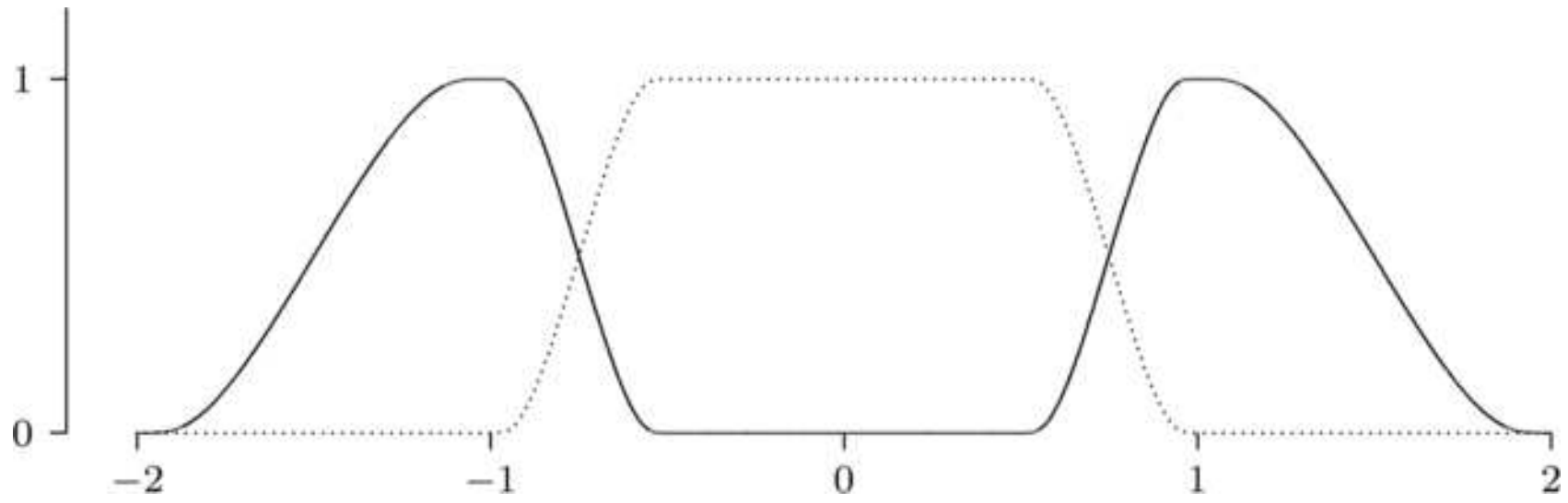

FIG. 1. *Typical graph of* $\phi$ *(dots) and* $b^2$ *(solid). Here* $B=2$.

2.1. *Littlewood–Paley decomposition.* Let $\phi$ be a $C^\infty$ function on $\mathbb{R}$, symmetric and decreasing on $\mathbb{R}^+$ supported in $|\xi| \le 1$, such that $1 \ge \phi(\xi) \ge 0$ and $\phi(\xi) = 1$ if $|\xi| \le \frac{1}{B}$. Let us define for an arbitrarily chosen $B > 1$ (see Figure 1):

$$b^2(\xi) = \phi\left(\frac{\xi}{B}\right) - \phi(\xi) \ge 0,$$

so that

$$\forall |\xi| \ge 1 \qquad \sum_{j \ge 0} b^2\left(\frac{\xi}{B^j}\right) = 1. \tag{4}$$

Remark that $b(\xi) \ne 0$ only if $\frac{1}{B} \le |\xi| \le B$. Let us now define the operator $\Lambda_j = \sum_{l \ge 0} b^2(\frac{l}{B^j}) L_l$ and the associated kernel

$$\Lambda_j(x,y) = \sum_{l \ge 0} b^2\left(\frac{l}{B^j}\right) L_l(\langle x,y\rangle) = \sum_{B^{j-1} < l < B^{j+1}} b^2\left(\frac{l}{B^j}\right) L_l(\langle x,y\rangle).$$

The following proposition is obvious.

PROPOSITION 1. *For every* $f \in \mathbb{L}^2$, $f = \lim_{J \to \infty} L_0(f) + \sum_{j=0}^{J} \Lambda_j(f)$, *where* $L_l(f) = \int_{\mathbb{S}^2} L_l(\langle x,y\rangle) f(y)\,dy$ *and* $\Lambda_j(f) = \int \Lambda_j(x,y) f(y)\,dy$. *Moreover, if* $M_j(x,y) = \sum_{l \ge 0} b(\frac{l}{B^j}) L_l(\langle x,y\rangle)$, *then*

$$\Lambda_j(x,y) = \int M_j(x,z) M_j(z,y)\,dz. \tag{5}$$

2.2. *Discretization and localization properties.* Let

$$\mathcal{K}_l = \bigoplus_{m=0}^{l} \mathcal{H}_m,$$

the space of the restrictions to $\mathbb{S}^2$ of the polynomials of degree less than $l$. The following quadrature formula is true: for all $l \in \mathbb{N}$ there exists a finite subset $\mathcal{X}_l \subset \mathbb{S}^2$ and positive real numbers $\lambda_\eta > 0$, indexed by $\eta \in \mathcal{X}_l$, such that

$$\forall f \in \mathcal{K}_l \qquad \int_{\mathbb{S}^2} f(x)\,dx = \sum_{\eta \in \mathcal{X}_l} \lambda_\eta f(\eta). \tag{6}$$

The operator $M_j$ defined in Proposition 1 is such that

$$z \mapsto M_j(x,z) \in \mathcal{K}_{[B^{j+1}]}, \qquad x \mapsto M_j(x,z) \in \mathcal{K}_{[B^{j+1}]},$$

so that

$$z \mapsto M_j(x,z)M_j(z,y) \in \mathcal{K}_{[2B^{j+1}]}$$

and, by the quadrature formula (6),

$$\Lambda_j(x,y) = \int M_j(x,z)M_j(z,y)\,dz = \sum_{\eta \in \mathcal{X}_{[2B^{j+1}]}} \lambda_\eta M_j(x,\eta)M_j(\eta,y).$$

This implies

$$\Lambda_j f(x) = \int \Lambda_j(x,y)f(y)\,dy = \int \sum_{\eta \in \mathcal{X}_{[2B^{j+1}]}} \lambda_\eta M_j(x,\eta)M_j(\eta,y)f(y)\,dy$$

$$= \sum_{\eta \in \mathcal{X}_{[2B^{j+1}]}} \sqrt{\lambda_\eta} M_j(x,\eta) \int \sqrt{\lambda_\eta} M_j(y,\eta)f(y)\,dy.$$

We denote

$$\mathcal{X}_{[2B^{j+1}]} = \mathcal{Z}_j, \qquad \psi_{j,\eta}(x) := \sqrt{\lambda_\eta} M_j(x,\eta)$$

and have

$$\frac{1}{c}B^{2j} \le \#\mathcal{Z}_j \le cB^{2j}$$

for some $c > 0$. We note $N_j = \sqrt{\#\mathcal{Z}_j}$. It holds, by Proposition 1,

$$f = L_0(f) + \sum_j \sum_{\eta \in \mathcal{Z}_j} \langle f, \psi_{j,\eta}\rangle_{\mathbb{L}^2} \psi_{j,\eta}.$$

The main result of Narcowich, Petrushev and Ward (2006) is the following localization property of the $\psi_{j,\eta}$, that are called "needlets": for any $k$ there exists a constant $c_k$ such that, for every $\xi \in \mathbb{S}^2$

$$|\psi_{j,\eta}(\xi)| \le \frac{c_k B^j}{(1 + B^j d(\eta,\xi))^k}, \tag{7}$$

where $d(\xi,\eta) = \arccos\langle \eta, \xi\rangle$ is the natural geodesic distance on the sphere. In other words, needlets are almost exponentially localized around any cubature point, which motivates their name (see Figure 2 in Section 8). Finally, notice that the construction in Narcowich, Petrushev and Ward (2006) is made with $B = 2$. We introduce here the free parameter $B > 1$, because in physical applications it may be useful in fine tuning the concentration in frequency.

**3. Isotropic fields on the sphere and the needlet expansion.** In this section we define the needlet expansion of an isotropic random field $T$ on $\mathbb{S}^2$. We say that $T$ is invariant by rotation (or isotropic) if

$$\forall \rho \in SO(3) \qquad \mathrm{E}[T(\rho x)T(\rho y)] = \mathrm{E}[T(x)T(y)].$$

This is equivalent to the fact that the covariance function of the process is of the form

$$\mathrm{E}[T(x)T(y)] = K(\langle x, y\rangle),$$

where $K$ is a bounded function defined on $[-1, +1]$.

Throughout this paper we make the following assumption.

ASSUMPTION 1. *$T$ is a centered Gaussian field, that is, mean square continuous and isotropic.*

Let us decompose $K$ on the basis of Legendre polynomials:

$$K = \sum_{l\ge 0} C_l L_l, \qquad C_l = \frac{8\pi^2}{2l+1}\int_{-1}^{1} K(u)L_l(u)\,du.$$

We write

$$T(x) = \sum_{l\ge 0} T_l(x),$$

where

$$T_l(x) = \int_{\mathbb{S}^2} T(y)L_l(\langle x, y\rangle)\,dx$$

($T_0 = 0$ as the field is assumed to be centered). It is immediate that

$$\mathrm{E}[T_l(x)T_k(y)] = \delta_{kl} C_l L_l(\langle x, y\rangle) \qquad \text{for every } x, y \in \mathbb{S}^2.$$

Actually all vectors in $H_l$ are eigenvectors of the Karhunen–Loève expansion of $K(\langle \cdot, \cdot\rangle)$. The previous projection can be realized explicitly as

$$T_l(x) = \sum_{m=-l}^{l} a_{lm} Y_{lm}(x), \tag{8}$$

where

$$a_{lm} = \int T(x) Y_{lm}(x)\,dx. \tag{9}$$

$(a_{lm})_{l,m}$ is a triangular array of complex uncorrelated [but for the condition $(-1)^m \overline{a}_{lm} = a_{l,-m}$] r.v.'s and $C_l$ is equal to the variance of $a_{lm}$.

For every integer $j$, let $\mathcal{Z}_j$ be the set of cubature points defined in the previous section. The points $\eta$ belonging to $\mathcal{Z}_j$ will be denoted $\xi_{jk}$, $k =$ $1,\dots,N_j$. Similarly we denote $\psi_{j,\eta}$ by $\psi_{j,k}$, and the needlet coefficient of a function $f$, $\langle f,\psi_{j,\eta}\rangle_{\mathbb{L}^2}$, by $\beta_{j,k}$. Hence the random needlet coefficients are

$$\beta_{j,k}:=\langle T,\psi_{j,k}\rangle_{\mathbb{L}^2}=\int_{\mathbb{S}^2}T(x)\psi_{j,k}(x)\,dx=\sqrt{\lambda_{j,k}}\sum_l b\left(\frac{l}{B^j}\right)T_l(\xi_{j,k}).$$

Actually

$$\begin{aligned}\int_{\mathbb{S}^2}T(x)\psi_{j,k}(x)\,dx&=\int_{\mathbb{S}^2}\sum_l T_l(x)\psi_{j,k}(x)\,dx\\&=\sqrt{\lambda_{j,k}}\sum_{l'}b\left(\frac{l'}{B^j}\right)\sum_l\int_{\mathbb{S}^2}T_l(x)L_{l'}(\langle x,\xi_{j,k}\rangle)\,dx\\&=\sqrt{\lambda_{j,k}}\sum_l b\left(\frac{l}{B^j}\right)\int_{\mathbb{S}^2}T_l(x)L_l(\langle x,\xi_{j,k}\rangle)\,dx\\&=\sqrt{\lambda_{j,k}}\sum_l b\left(\frac{l}{B^j}\right)T_l(\xi_{j,k}),\end{aligned}$$

in view of the reproducing properties of the projection kernel. Hence

$$\mathrm{E}[\beta_{j,k}\beta_{j,k'}]=\sqrt{\lambda_{j,k}\lambda_{j,k'}}\sum_l b^2\left(\frac{l}{B^j}\right)K_l(\langle\xi_{j,k},\xi_{j,k'}\rangle)$$

and

$$\begin{aligned}\mathrm{Cor}(\beta_{j,k},\beta_{j,k'})&:=\frac{\mathrm{E}[\beta_{j,k}\beta_{j,k'}]}{\sqrt{\mathrm{E}[\beta_{j,k}^2]\,\mathrm{E}[\beta_{j,k'}^2]}}\\&=\frac{\sqrt{\lambda_{j,k}\lambda_{j,k'}}}{\sqrt{\lambda_{j,k}\lambda_{j,k'}}}\frac{\sum_l b^2(l/B^j)K_l(\langle\xi_{j,k},\xi_{j,k'}\rangle)}{\sum_l b^2(l/B^j)C_lL_l(1)}\\&=\frac{\sum_l b^2(l/B^j)K_l(\langle\xi_{j,k},\xi_{j,k'}\rangle)}{\sum_l b^2(l/B^j)C_lL_l(1)}.\end{aligned}$$

We shall need to assume some regularity conditions on the asymptotic behavior of angular power spectrum $C_l$.

Assumption 2. There exist $M>0$, $\alpha>2$ and a sequence of functions $(g_j)_j$ such that

$$C_l=l^{-\alpha}g_j\left(\frac{l}{B^j}\right)$$

for every $l$ such that $B^{j-1}<l<B^{j+1}$, and positive numbers $c_1,c_2,k_r,r=0,\dots,M$, such that

$$c_1l^{-\alpha}\le C_l\le c_2l^{-\alpha},\qquad\sup_j\sup_{B^{-1}\le u\le B}|g_j^{(r)}(u)|\le k_r,$$

where $g_j^{(r)}(u) := \frac{d^r}{du^r} g_j(u)$ (a uniform bounded differentiability condition).

REMARK 2. Assumption 2 is a regularity condition on the asymptotics of the angular power spectrum which is trivially satisfied, for instance, if $C_l = l^{-\alpha}$. Note that the sequence $(g_j)_j$ belongs uniformly to the Sobolev space $W^{M,\infty}$.

The following result is the basic localization inequality which plays a crucial role for the arguments below.

LEMMA 3.

$$|\operatorname{Cor}(\beta_{j,k}, \beta_{j,k'})| \le \frac{C_M}{(1 + B^j d(\xi_{j,k}, \xi_{j,k'}))^M}, \tag{10}$$

*where, as hinted above,* $d(\xi_{j,k}, \xi_{j,k'}) = \arccos(\langle \xi_{j,k}, \xi_{j,k'} \rangle)$.

PROOF. Observe first that, as we assumed that $c_1 l^{-\alpha} \le C_l \le c_2 l^{-\alpha}$,

$$c_1 B^{(2-\alpha)j} \le \sum_l b^2\left(\frac{l}{B^j}\right) C_l L_l(1) \le c_2 B^{(2-\alpha)j}. \tag{11}$$

We recall the following bound for type II polynomials which is derived in Narcowich, Petrushev and Ward (2006), Theorem 2.6:

$$\left|\sum_l \phi_j\left(\frac{l}{B^j}\right) L_l(\langle x, y\rangle)\right| \le \frac{c_M B^{2j}}{(1 + B\, d(x,y))^M},$$

where $c_M$ only depends on $\sup_{j\ge 1, k\le M} \|\phi_j^{(k)}\|_1$. Whence, using this for $\phi_j(x) = b^2(x) x^{-\alpha} g_j(x)$ and (11),

$$\begin{aligned}
&\left|\frac{\sum_l b^2(l/B^j) K_l(\langle \xi_{j,k}, \xi_{j,k'}\rangle)}{\sum_l b^2(l/B^j) C_l L_l(1)}\right| \\
&\quad = \left|\frac{\sum_l b^2(l/B^j) C_l L_l(\langle \xi_{j,k}, \xi_{j,k'}\rangle)}{\sum_l b^2(l/B^j) C_l L_l(1)}\right| \\
&\quad = \left|\frac{\sum_l b^2(l/B^j) l^{-\alpha} g_j(l/B^j) L_l(\langle \xi_{j,k}, \xi_{j,k'}\rangle)}{\sum_l b^2(l/B^j) C_l L_l(1)}\right| \\
&\quad = \left|\frac{\sum_l b^2(l/B^j)(l/B^j)^{-\alpha} g_j(l/B^j) L_l(\langle \xi_{j,k}, \xi_{j,k'}\rangle)}{B^{j\alpha}\sum_l b^2(l/B^j) C_l L_l(1)}\right| \\
&\quad \le \frac{C_M}{(1 + B^j d(\xi_{j,k}, \xi_{j,k'}))^M}.
\end{aligned}$$

□

REMARK 4. As mentioned in the Introduction, the previous lemma highlights a peculiar feature: the needlet coefficients at any finite distance are asymptotically uncorrelated. This property is at the heart of our results below.

**4. The central limit theorem for polynomial functionals of needlet coefficients.** In the sequel, we make an extensive use of *diagrams*, which are mnemonic devices for computation of moments and cumulants of polynomials of Gaussian random variables. We adopt standard terms and notation (edges, nodes, connected components, ...), and we refer to Surgailis (2003), Marinucci (2006) and Baldi et al. (2008), Section 5, for definitions and background results. See also Nualart and Peccati (2005) for a more recent point of view.

In the arguments to follow, we focus on polynomial functionals of the (normalized) wavelets coefficients, of the form

$$h_{u,N_j} := \frac{1}{N_j}\sum_{k=1}^{N_j^2}\sum_{q=1}^{Q} w_{uq} H_q(\widehat{\beta}_{j,k}), \qquad \widehat{\beta}_{j,k} := \frac{\beta_{j,k}}{\sqrt{\mathrm{E}\,\beta_{j,k}^2}}, \tag{12}$$

where $u = 1, 2, \ldots, U$. Recall that $N_j = \sqrt{\#\mathcal{Z}_j}$. Here $w_{uq}$ are real scalars and $H_q$ denotes the $q$th Hermite polynomial. As Hermite polynomials are an algebraic basis, every polynomial in the variables $\widehat{\beta}_{j,k}$ is of the form (12); we start from a general characterization of the behavior of the sequences $h_{u,N_j}$. First we define the covariance matrix $\Omega_j$, with elements

$$\Omega_j := \{\mathrm{E}[h_{u,N_j} h_{v,N_j}]\}_{u,v=1,\ldots,U}. \tag{13}$$

Throughout the sequel we assume the following regularity condition:

ASSUMPTION 3. There exists $j_0$ such that for $j \geq j_0$ the covariance matrix $\Omega_j$ is invertible.

Assumption 3 is a nondegeneracy condition on the asymptotics of the statistics of interest. Consider for instance the scalar case $U = 1$. From the diagram formula, it is immediate to obtain

$$\begin{aligned}\mathrm{E}[h_{u,N_j}^2] &= \frac{1}{N_j^2}\sum_{q=1}^{Q} w_q^2 \operatorname{Var}\left(\sum_k H_q(\widehat{\beta}_{j,k})\right)\\ &= \frac{1}{N_j^2}\sum_{q=1}^{Q} q!\, w_{uq}^2 \sum_{k,k'=1}^{N_j^2}\left[\frac{\mathrm{E}[\beta_{j,k}\beta_{j,k'}]}{\sqrt{\mathrm{E}[\beta_{j,k}^2]\,\mathrm{E}[\beta_{j,k'}^2]}}\right]^q.\end{aligned}$$

The previous condition merely states that our nonlinear statistics have a nondegenerate asymptotic variance. Ruling aside pathological cases, it should be noted that the previous assumption basically requires $w_{up}^2 > 0$ for some $p$. In the multivariate case $U > 1$ we also require that the polynomials $h_{u,N_j}$ and $h_{v,N_j}$ are linearly independent. It is to be noted, however, that the assumption fails for a polynomial of order 1.

THEOREM 5. *Under Assumptions 1, 2 and 3, as* $N_j \to \infty$

$$\Omega_j^{-1/2}(h_{1N_j}, \ldots, h_{UN_j})' \underset{j\to\infty}{\overset{\mathcal{D}}{\to}} N(0, I_U),$$

*where* $I_U$ *denotes the identity matrix of dimension* $U$.

PROOF. We note first that the multivariate result follows immediately from the case $U = 1$, as by the Cramér–Wald device it is enough to focus on sequences of the form

$$h_{N_j} := \sum_{u=1}^{U} \lambda_u h_{u,N_j}$$

and

$$\widehat{h}_{N_j} := \frac{\sum_{u=1}^{U} \lambda_u h_{u,N_j}}{\sqrt{\sum_{u,v=1}^{U} \lambda_u \lambda_v \operatorname{E}[h_{u,N_j} h_{v,N_j}]}} = \frac{\sum_{u=1}^{U} \lambda_u h_{u,N_j}}{\sqrt{\operatorname{E}[(\sum_{u=1}^{U} \lambda_u h_{u,N_j})^2]}}.$$

However, it is clear that, for any choice of real numbers $\lambda_1, \ldots, \lambda_U$,

$$h_{N_j} = \frac{1}{N_j} \sum_{k=1}^{N_j^2} \sum_{q=1}^{Q} \sum_{u=1}^{U} \lambda_u w_{uq} H_q(\beta_{j,k}) = \frac{1}{N_j} \sum_{k=1}^{N_j^2} \sum_{q=1}^{Q} \widetilde{w}_q H_q(\beta_{j,k}),$$

where $\widetilde{w}_q := \sum_{u=1}^{U} \lambda_u w_{uq}$. It is obvious that $\operatorname{E}[h_{u,N_j}] = 0$. Hence to complete the argument it is sufficient to prove that, as $N_j \to \infty$,

$$\lim_{N\to\infty} \operatorname{E}\left[\left(\frac{h_{N_j}}{\sqrt{\operatorname{Var}(h_{N_j})}}\right)^p\right] = \begin{cases} (p-1)!!, & \text{for } p = 2, 4, \ldots, \\ 0, & \text{otherwise.} \end{cases}$$

We must show that, as $N_j \to \infty$ and for all $p \geq 3$,

$$\begin{aligned}
&\operatorname{Cum}\left(\frac{1}{N_j} \sum_{k_1=1}^{N_j^2} \sum_{q=1}^{Q} \widetilde{w}_q H_q(\widehat{\beta}_{j,k_1}), \ldots, \frac{1}{N_j} \sum_{k_p=1}^{N_j^2} \sum_{q=1}^{Q} \widetilde{w}_q H_q(\widehat{\beta}_{j,k_p})\right) \\
&\quad = \sum_{q_1,\ldots,q_p}^{Q} \widetilde{w}_{q_1} \cdots \widetilde{w}_{q_p} \operatorname{Cum}\left(\frac{1}{N_j} \sum_{k_1=1}^{N_j^2} H_{q_1}(\widehat{\beta}_{j,k_1}), \ldots, \frac{1}{N_j} \sum_{k_p=1}^{N_j^2} H_{q_p}(\widehat{\beta}_{j,k_p})\right) \\
&\quad = \frac{1}{N_j^p} \sum_{q_1,\ldots,q_p}^{Q} \widetilde{w}_{q_1} \cdots \widetilde{w}_{q_p} \sum_{G \in \mathcal{V}_c(q_1,\ldots,q_p)} \sum_{k_1,\ldots,k_p=1}^{N_j^2} \prod_{1\leq u<v\leq p} \gamma_{k_u k_v}^{\eta_{uv}(G)}
\end{aligned}$$

$$\leq \frac{c}{N_j^p} \sup_{q_1,\dots,q_p} \sum_{G\in\mathcal{V}_c(q_1,\dots,q_p)} \sum_{k_1,\dots,k_p=1}^{N_j^2} \prod_{1\leq u<v\leq p} |\gamma_{k_v k_u}|^{\eta_{uv}(G)} \to 0,$$

where $\eta_{uv}(G)$ counts the number of edges between node $u$ and node $v$. It is then clearly enough if we can prove that

$$\begin{aligned}\sum_{k_1,\dots,k_p=1}^{N_j^2} \prod_{1\leq u<v\leq p} |\gamma_{k_u k_v}|^{\eta_{uv}(G)} \\ &= O(N_j^{2[(p-1)/2]}) \\ &= \begin{cases} O(N_j^{p-1}), & \text{for } p=3,5,7,\dots, \\ O(N_j^{p-2}), & \text{for } p=2,4,6,\dots. \end{cases}\end{aligned}$$

Now write

$$\chi_{q_1\cdots q_p}(G) := \sum_{k_1,\dots,k_p=1}^{N_j^2} \prod_{1\leq u<v\leq p} |\gamma_{k_u k_v}|^{\eta_{uv}(G)}.$$

Note that each of the covariances is bounded by 1, so that $\chi_{q_1\cdots q_p}(G)$ is a nonincreasing function of $\eta_{uv}(G)$, $u,v=1,\dots,p$. We modify iteratively the elements $\eta_{uv}(G)$ by picking $(u,v)$ at random, and then decreasing $\eta_{uv}(G)$ by 1; in graphical terms, this can be viewed as taking a new graph $G_1$ where an edge between $u$ and $v$ has been deleted ($G_1$ need no longer be connected). We repeat this procedure until (in a finite number of steps, $T$, say), we obtain a graph, $G_T$, such that the following circumstances are met:

(a) There are no isolated nodes.

(b) There exists at least a path covering three nodes.

(c) The connected components do not allow loops.

It is simple to see that we can reach $G_T$ in a finite number of steps by the following algorithm:

(1) We keep lowering $\eta_{uv}$ until we get to the point where the next step would necessarily violate condition (a)

(2) If condition (b) is met, we stop our procedure.

(3) If condition (b) fails, it means we have only components with two nodes and it is sufficient to raise by a unity any of the $\eta_{uv}$ (i.e., to introduce an edge between two components).

It is clear that there are at most $[\frac{p-1}{2}]$ such components. For brevity we assume that there are no paths with more than three nodes, the argument in the remaining cases being entirely analogous. We partition the nodes $u=1,\dots,p$ into subsets $I_1$ and $I_2$ according to the following rule. All nodes that belong to more than one edge belong to $I_1$; then for components with only two nodes we put the one whose index is smaller again into $I_1$. All the remaining others are put into $I_2$. It is simple to check that the cardinality of $I_1$ equals the number of unconnected components in $G_T$ and hence is smaller than $\frac{p-1}{2}$. Since $|\gamma_{k_u k_v}| \le 1$, we have

$$\prod_{1\le u<v\le p} |\gamma_{k_u k_v}|^{\eta_{uv}(G_T)} = \prod_{\substack{u\in I_1\cup I_2\\ u<v\le p}} |\gamma_{k_u k_v}|^{\eta_{uv}(G_T)} \le \prod_{\substack{u\in I_2\\ u<v\le p}} |\gamma_{k_u k_v}|^{\eta_{uv}(G_T)} \le \prod_{\substack{u\in I_2\\ u<v\le p}} |\gamma_{k_u k_v}|.$$

Note that by construction, $k_u$ appears exactly once in the covariances whenever $u \in I_2$; hence we obtain

$$\sum_{k_1,\dots,k_p=1}^{N_j^2} \prod_{\substack{u\in I_2\\ u<v\le p}} |\gamma_{k_u k_v}| = \sum_{k_u,u\in I_1} \sum_{k_u,u\in I_2} \prod_{\substack{u\in I_2\\ u<v\le p}} |\gamma_{k_u k_v}| = \sum_{k_u,u\in I_1} \prod_{\substack{u\in I_2\\ u<v\le p}} \left( \sum_{k_u,u\in I_2} |\gamma_{k_u k_v}| \right).$$

Thus we obtain, using (10) and the following Lemma 6,

$$\sum_{k_u,u\in I_1} \prod_{\substack{u\in I_2\\ u<v\le p}} \left( \sum_{k_u,u\in I_2} |\gamma_{k_u k_v}| \right) \le \sum_{k_u,u\in I_1} \left[ \sum_{k'\in \mathcal{Z}_j} \frac{C_M}{(1+B^j d(\xi_{j,k},\xi_{j,k'}))^M} \right]^{p^2} \le \sum_{k_u,u\in I_1} C = O(N_j^{2[(p-1)/2]}).$$

□

LEMMA 6. *If $M \ge 3$, there exists a constant $C'_M$ such that*

$$\sum_{k'=1}^{N_j^2} \frac{1}{(1+B^j d(\xi_{j,k},\xi_{j,k'}))^M} \le C'_M.$$

PROOF. It is proved in Narcowich, Petrushev and Ward (2006) that to get cubature points for polynomials of degree less than $L$, it is enough to take a maximal $\epsilon$-mesh on the sphere with $\epsilon \sim \frac{1}{L}$ [i.e., a set $\{x_1, \dots, x_K\}$ with $d(x_i, x_j) > \epsilon$ for $x_i \neq x_j$ and $K$ maximal]. Using a simple covering argument, we have $\frac{2\pi}{|B(\epsilon)|} \leq K \leq \frac{2\pi}{|B(\epsilon/2)|}$ where $|B(\epsilon)|$ is the volume of (any) ball of radius $\epsilon$, and $|B(\epsilon)| \sim \frac{1}{\epsilon^2}$, so $K \sim L^2$.

Let $L \sim B^j$ and the corresponding mesh defining $\mathcal{Z}_j$; as balls are disjoint and $d(\xi_{jk}, x) \leq 2d(\xi_{jk}, \xi_{jk'})$ by the triangular inequality, we obtain

$$
\begin{aligned}
&\sum_{k'=1}^{N_j^2} \frac{1}{(1 + B^j d(\xi_{j,k}, \xi_{j,k'}))^M} \\
&\quad = \frac{1}{|B(1/(2L))|} \sum_{k' \in \mathcal{Z}_j} \int_{B(\xi_{jk'}, 1/(2L))} \frac{dx}{(1 + Ld(\xi_{j,k}, \xi_{j,k'}))^M} \\
&\quad \leq CL^2 \sum_{k' \in \mathcal{Z}_j} \int_{B(\xi_{jk'}, 1/(2L))} \frac{dx}{(1 + Ld(\xi_{j,k}, x))^M} \\
&\quad \leq CL^2 \int_{\mathbb{S}^2} \frac{dx}{(1 + L \arccos(\langle \xi_{j,k}, x \rangle))^M} = 2C\pi L^2 \int_0^\pi \frac{\sin\theta \, d\theta}{(1 + L\theta)^M} \\
&\quad \leq 2C\pi L^2 \int_0^\pi \frac{\theta \, d\theta}{(1 + L\theta)^M} \leq 2C\pi L^2 \left( \int_0^{1/L} \theta \, d\theta + \frac{1}{L^M} \int_{1/L}^\infty \theta^{1-M} \, d\theta \right) \\
&\quad = 2C\pi L^2 \left( \frac{1}{2L^2} + \frac{1}{L^M} \frac{1}{(M-2)L^{2-M}} \right) \leq 2C\pi. \qquad \square
\end{aligned}
$$

**5. The functional central limit theorem.** We are now ready to introduce the following continuous-time vector process:

$$
W_J(r) := \frac{1}{\sqrt{J}} \sum_{j=2,4,\dots}^{[Jr]} \Omega_j^{-1/2} (h_{1,N_j}, \dots, h_{U,N_j})', \qquad 0 \leq r \leq 1,
$$

where $\Omega_j$ was defined in (13). Here $U \geq 1$ is a fixed integer.

THEOREM 7. *Under Assumptions 1, 2 and 3, as* $N_j \to \infty$

$$
W_J \quad \Longrightarrow \quad X,
$$

*where* $X$ *denotes the* $U$*-dimensional standard Brownian motion and* $\Longrightarrow$ *denotes weak convergence in the Skorohod space* $D([0,1]^U)$.

PROOF. We note first that the multivariate result follows from the case $U = 1$, as remarked above. It is well known that in order to prove weak convergence we have to establish convergence of the finite-dimensional distributions and tightness. By the Cramér–Wald device, to establish the former it is enough to focus on sequences of the form

$$\frac{1}{\sqrt{J}}\sum_{j=2,4,\ldots}^{[Jr]} h_{N_j} := \frac{1}{\sqrt{J}}\sum_{j=2,4,\ldots}^{[Jr]}\sum_{u=1}^{U}\lambda_u h_{u,N_j},$$

$$\frac{1}{\sqrt{J}}\sum_{j=2,4,\ldots}^{[Jr]} \widehat{h}_{N_j} := \frac{\sum_{j=2,4,\ldots}^{[Jr]}\sum_{u=1}^{U}\lambda_u h_{u,N_j}}{\sqrt{\sum_{j=2,4,\ldots}^{[Jr]}\sum_{u,v=1}^{U}\lambda_u\lambda_v \operatorname{E} h_{u,N_j}h_{v,N_j}}}$$

$$= \frac{\sum_{j=2,4,\ldots}^{[Jr]}\sum_{u=1}^{U}\lambda_u h_{u,N_j}}{\sqrt{\operatorname{E}(\sum_{j=2,4,\ldots}^{[Jr]}\sum_{u=1}^{U}\lambda_u h_{u,N_j})^2}}.$$

However, it is clear that for any choice of real numbers $\lambda_1,\ldots,\lambda_U$,

$$\frac{1}{\sqrt{J}}\sum_{j=2,4,\ldots}^{[Jr]}\left\{\frac{1}{N_j}\sum_{k=1}^{N_j^2}\sum_{q=1}^{Q}\sum_{u=1}^{U}\lambda_u w_{uq}H_q(\beta_{j,k})\right\}$$

$$= \frac{1}{\sqrt{J}}\sum_{j=2,4,\ldots}^{[Jr]}\left\{\frac{1}{N_j}\sum_{k=1}^{N_j^2}\sum_{q=1}^{Q}\widetilde{w}_q H_q(\beta_{j,k})\right\},$$

where, as before, $\widetilde{w}_q := \sum_{u=1}^{U}\lambda_u w_{uq}$. On the other hand, a necessary and sufficient condition for tightness of vector processes is tightness for the component processes. Without any loss of generality, we can hence focus on the univariate case $U=1$. We first consider convergence of the finite-dimensional distributions. It is straightforward to see that $\widehat{h}_{N_j}, \widehat{h}_{N_{j'}}$ are independent whenever $|j-j'|\geq 2$. As the process $W_J(r)$ is a partial sum of independent elements, convergence of the finite-dimensional distributions follows from the Lyapunov condition

$$\lim_{J\to\infty}\frac{\sum_{j=2,4,\ldots}^{[Jr]}\operatorname{E}[\widehat{h}_{N_j}^4]}{J^2r^2} = 0. \tag{14}$$

We have

$$\operatorname{E}[\widehat{h}_{N_j}^4] = \operatorname{E}\left[\frac{1}{N_j}\sum_k\sum_{q=1}^{Q}\widetilde{w}_q H_q(\widehat{\beta}_{j,k})\right]^4$$

$$= \sum_{q_1q_2q_3q_3}^{Q}\widetilde{w}_{q_1}\widetilde{w}_{q_2}\widetilde{w}_{q_3}\widetilde{w}_{q_4}$$

$$\times\left[\sum_{G\in\mathcal{V}(q_1,\ldots,q_4)}\frac{1}{N_j^4}\sum_{k_1k_2k_3k_4}\prod_{1\leq u<v\leq 4}|\operatorname{E}[\widehat{\beta}_{j,k_u}\widehat{\beta}_{j,k_v}]|^{\eta_{uv}(G)}\right]$$

$$\begin{aligned}
&\le \frac{c}{N_j^4} \sum_{k_1 k_2 k_3 k_4} \{|\operatorname{E}[\widehat{\beta}_{j,k_1}\widehat{\beta}_{j,k_2}]\operatorname{E}[\widehat{\beta}_{j,k_3}\widehat{\beta}_{j,k_4}]| \\
&\qquad + |\operatorname{E}[\widehat{\beta}_{j,k_1}\widehat{\beta}_{j,k_3}]\operatorname{E}[\widehat{\beta}_{j,k_2}\widehat{\beta}_{j,k_4}]| \\
&\qquad + |\operatorname{E}[\widehat{\beta}_{j,k_1}\widehat{\beta}_{j,k_4}]\operatorname{E}[\widehat{\beta}_{j,k_2}\widehat{\beta}_{j,k_4}]|\} \\
&\le c\left[\frac{1}{N_j^2}\sum_{k_1 k_2} |\operatorname{E}[\widehat{\beta}_{j,k_1}\widehat{\beta}_{j,k_2}]|\right]^2 = O(1),
\end{aligned}$$

uniformly over $j$, in view of Lemma 6. Equation (14) then follows easily from

$$\frac{\sum_{j=2,4,\dots}^{[Jr]} \operatorname{E}\widehat{h}_{N_j}^4}{J^2 r^2} \le c\frac{[Jr]}{[Jr]^2} \underset{J\to\infty}{\to} 0.$$

Likewise, by a well-known result, tightness follows from

$$\begin{aligned}
&\operatorname{E}[|W_J(r_1) - W_J(r)|^2 |W_J(r_2) - W_J(r)|^2] \\
&\quad = \frac{1}{J^2}\operatorname{E}\left[\left(\sum_{j=[Jr]+1}^{[Jr_2]} \widehat{h}_{N_j}\right)^2\right]\operatorname{E}\left[\left(\sum_{j=[Jr_1]+1}^{[Jr]} \widehat{h}_{N_j}\right)^2\right] \\
&\quad \le \frac{C}{J^2}([Jr_2]-[Jr])([Jr]-[Jr_1]) \le 4C(r_2-r_1)^2
\end{aligned}$$

for all $r_1 \le r \le r_2$, again in view of Lemma 6. □

**6. Statistical applications.** In this section, we use the previous results to derive goodness-of-fit for spherical random fields. In particular, we take

$$\begin{aligned}
h_{1N_j} &= \frac{1}{N_j}\sum_{k=1}^{N_j^2} H_2(\widehat{\beta}_{j,k}) = \frac{1}{N_j}\sum_{k=1}^{N_j^2}\{\widehat{\beta}_{j,k}^2 - 1\}, \\
h_{2N_j} &= \frac{1}{N_j}\sum_{k=1}^{N_j^2}\{H_3(\widehat{\beta}_{j,k}) + 3H_1(\widehat{\beta}_{j,k})\} = \frac{1}{N_j}\sum_{k=1}^{N_j^2}\widehat{\beta}_{j,k}^3, \\
h_{3N_j} &= \frac{1}{N_j}\sum_{k=1}^{N_j^2}\{H_4(\widehat{\beta}_{j,k}) + 6H_2(\widehat{\beta}_{j,k})\} = \frac{1}{N_j}\sum_{k=1}^{N_j^2}\{\widehat{\beta}_{j,k}^4 - 3\}.
\end{aligned}$$

It is natural to view $h_{1N_j}$ as a goodness-of-fit statistic on the angular power spectrum $\{C_l\}$. More precisely, a typical question arising in applications is to check the validity of a physical model (e.g., specific values of parameters) by means of a comparison between the expected and observed angular power spectrum. In this framework, this goal can be accomplished as follows: recall

that

$$\frac{1}{N_j}\sum_{k=1}^{N_j^2}\{\widehat{\beta}_{j,k}^2-1\}=\frac{1}{N_j}\sum_{k=1}^{N_j^2}\frac{\beta_{j,k}^2-\mathrm{E}[\beta_{j,k}^2]}{\mathrm{E}[\beta_{j,k}^2]},$$

where

$$\mathrm{E}[\beta_{j,k}^2]=\frac{1}{N_j^2}\sum_{B^{j-1}\le l\le B^{j+1}}b^2\biggl(\frac{l}{B^j}\biggr)C_l\frac{2l+1}{4\pi}.$$

Then it is clear that $h_{1N_j}$ provides a measure of discrepancy between the expected and observed values of an averaged power spectrum. In order to construct feasible statistical procedures with an asymptotic justification to investigate this and other hypotheses of interest, we define

$$W_J(r)=\frac{1}{\sqrt{J}}\sum_{j=B,B^2,\dots}^{[Jr]}\Omega_j^{-1/2}(h_{1,N_j},h_{2,N_j},h_{3,N_j})',\qquad 0\le r\le 1,$$

where as before

$$\Omega_j=\{\mathrm{E}[h_{u,N_j}h_{v,N_j}]\}_{u,v=1,\dots,U}$$

and

$$\mathrm{E}[h_{1,N_j}^2]=\mathrm{Var}(h_{1,N_j})=\frac{1}{N_j^2}\mathrm{Var}\biggl(\sum_k H_2(\widehat{\beta}_{j,k})\biggr)=\frac{2}{N_j^2}\sum_{k_1k_2}(\mathrm{E}[\widehat{\beta}_{j,k_1}\widehat{\beta}_{j,k_2}])^2,$$

$$\begin{aligned}\mathrm{E}[h_{2,N_j}^2]&=\mathrm{Var}(h_{2,N_j})=\frac{1}{N_j^2}\mathrm{Var}\biggl(\sum_k\{H_3(\widehat{\beta}_{j,k})+3H_1(\widehat{\beta}_{j,k})\}\biggr)\\&=\frac{6}{N_j^2}\sum_{k_1k_2}(\mathrm{E}[\widehat{\beta}_{j,k_1}\widehat{\beta}_{j,k_2}])^3+\frac{9}{N_j^2}\sum_{k_1k_2}\mathrm{E}[\widehat{\beta}_{j,k_1}\widehat{\beta}_{j,k_2}],\\\mathrm{E}[h_{3,N_j}^2]&=\mathrm{Var}(h_{3,N_j})=\frac{1}{N_j^2}\mathrm{Var}\biggl(\sum_k\{H_4(\widehat{\beta}_{j,k})+6H_2(\widehat{\beta}_{j,k})\}\biggr)\\&=\frac{24}{N_j^2}\sum_{k_1k_2}(\mathrm{E}[\widehat{\beta}_{j,k_1}\widehat{\beta}_{j,k_2}])^4+\frac{72}{N_j^2}\sum_{k_1k_2}(\mathrm{E}[\widehat{\beta}_{j,k_1}\widehat{\beta}_{j,k_2}])^2.\end{aligned}$$

Also

$$\begin{aligned}\mathrm{E}[h_{1,N_j}h_{2,N_j}]&=\frac{1}{N_j^2}\sum_{k_1k_2}\mathrm{E}[\{H_3(\widehat{\beta}_{j,k_1})+3H_1(\widehat{\beta}_{j,k_1})\}H_2(\widehat{\beta}_{j,k_2})]=0,\\\mathrm{E}[h_{1,N_j}h_{3,N_j}]&=\frac{1}{N_j^2}\sum_{k_1k_2}\mathrm{E}[\{H_4(\widehat{\beta}_{j,k_1})+6H_2(\widehat{\beta}_{j,k_1})\}H_2(\widehat{\beta}_{j,k_2})]\\&=\frac{12}{N_j^2}\sum_{k_1k_2}(\mathrm{E}[\widehat{\beta}_{j,k_1}\widehat{\beta}_{j,k_2}])^2\end{aligned}$$

and

$$\mathrm{E}[h_{2,N_j}h_{3,N_j}] = \frac{1}{N_j^2}\sum_{k_1k_2}\mathrm{E}[\{H_3(\widehat{\beta}_{j,k_1})+3H_1(\widehat{\beta}_{j,k_1})\}\{H_4(\widehat{\beta}_{j,k_1})+6H_2(\widehat{\beta}_{j,k_1})\}]=0.$$

As from Theorem 7, as $J\to\infty$ $W_J$ converges in $D([0,1]^3)$ to a three-dimensional standard Brownian motion $X$, by focusing on the first row of $W_J$, we obtain for instance the Kolmogorov–Smirnov type test, that is,

$$\lim_{J\to\infty} P\Big(\sup_{0\le r\le 1}|W_{1,J}(r)|>t\Big)=P\Big(\sup_{0\le r\le 1}|X_1(r)|>t\Big),$$

$X_1$ denoting the first component of $X$. The derivation of threshold values for $t$ is then standard. Similarly, it is possible to construct tests for Gaussianity and isotropy based on the skewness and kurtosis statistics $h_{2N_j}$ and $h_{3N_j}$, respectively. The numerical implementation of these procedures on CMB data is currently underway.

**7. Missing observations.** As mentioned in the Introduction, we expect needlets to be extremely robust in the presence of partially observed spherical random fields, due to their excellent localization properties in real space. This result can be formalized as follows: we assume we observe $\widetilde{T}(\xi)=T(\xi)+V(\xi)$, where $V(\xi)$ is a noise field that need not be independent from $T(\xi)$; indeed the most relevant case is $V(\xi)=-T(\xi)1_{\{\xi\in G\}}$, $G\subset\mathbb{S}^2$ denoting the unobserved subset of the sphere. This situation arises when the field is not observed (and hence its value is set to zero) for some locations in the sky. This is the situation with CMB data in the so-called *galactic cut* region, where CMB is dominated by the Milky Way emissions. We note $N_\varepsilon(\xi_{j,k}):=\{\xi\in\mathbb{S}^2:d(\xi,\xi_{j,k})\le\varepsilon\}$ the neighborhood of radius $\varepsilon$ around the cubature point $\xi_{j,k}$, $d$ denoting as usual the angular distance. We write

$$\widetilde{\beta}_{jk}=\int_{\mathbb{S}^2}\widetilde{T}(\xi)\psi_{jk}(\xi)\,d\xi$$

for the wavelets coefficients of $\widetilde{T}$. The following result highlights the robustness property of needlets.

PROPOSITION 8. *Let $\xi_{j,k}$ be a cubature point such that $V(\xi)=0$ on $N_\varepsilon(\xi_{j,k})$ and assume that*

$$\sup_{\xi\in\mathbb{S}^2}\mathrm{E}[V(\xi)^2]=:V^*<\infty. \tag{15}$$

*Then, for every $M\in\mathbb{N}$,*

$$\|\widetilde{\beta}_{jk}-\beta_{jk}\|_2:=\sqrt{\mathrm{E}[(\widetilde{\beta}_{jk}-\beta_{jk})^2]}\le\frac{C_M4\pi\sqrt{2V^*}B^j}{(1+B^j\varepsilon)^M}.$$

PROOF. We have, by the localization property (7),

$$
\begin{aligned}
(\widetilde{\beta}_{jk}-\beta_{jk})^2 &\le \left(\int_{\mathbb{S}^2\setminus N_\varepsilon(\xi_{j,k})} V(\xi)\psi_{jk}(\xi)\,d\xi\right)^2\\
&\le \left(\frac{C_M B^j}{(1+B^j\varepsilon)^M}\right)^2\left(\int_{\mathbb{S}^2}|V(\xi)|\,d\xi\right)^2\\
&\le 4\pi\left(\frac{C_M B^j}{(1+B^j\varepsilon)^M}\right)^2\left(\int_{\mathbb{S}^2}|V(\xi)|^2\,d\xi\right).
\end{aligned}
$$

Therefore

$$
\mathrm{E}[(\widetilde{\beta}_{jk}-\beta_{jk})^2] \le 4\pi\int_{\mathbb{S}^2}\mathrm{E}[V(\xi)^2]\,d\xi \le 16\pi^2 V^*\left(\frac{C_M B^j}{(1+B^j\varepsilon)^M}\right). \qquad \square
$$

REMARK 9. Remark that in Proposition 8 $V$ is not assumed to be isotropic. Thus in the case of gaps (15) is obviously satisfied, as $\mathrm{E}[V(\xi)^2]\le\mathrm{E}[T(\xi)^2]$ for every $\xi\in S^2$. It is also interesting to stress that, in view of (11),

$$
\sqrt{\frac{\mathrm{E}[(\widetilde{\beta}_{jk}-\beta_{jk})^2]}{\mathrm{E}[\beta_{jk}^2]}} \le c_1 B^{j(\alpha-1)}\frac{C_M B^j}{(1+B^j\varepsilon)^M}2\pi\sqrt{2V^*} \le C'_M B^{2j(\alpha-M)}
$$

with $C'_M=\frac{1}{\varepsilon}c_1 C_M 2\pi\sqrt{2V^*}$. For $M$ large enough, it is not difficult to show that, up to different normalizing constants, the limit results in Sections 4, 5 and 6 are not affected asymptotically by the presence of sky cuts. Although this result must be taken with a good deal of common sense when working with finite-resolution experiments, we view this property as a very strong rationale to motivate the use of needlets in cosmology and astrophysics.

## 8. Numerical implementation.

In this section we address some practical issues concerning the implementation of needlets on real data. In particular we consider data on the Cosmic Microwave Background radiation, as provided by the NASA experiment WMAP. It is not difficult to devise some kernel construction that fulfills the conditions highlighted in Section 2. As in Marinucci et al. (2008), we suggest the following algorithm [cf. Guilloux, Fay and Cardoso (2007) for alternative suggestions].

In order to construct the function $\varphi$ of Section 2.1 one just defines $f(t)=\exp(-\frac{1}{1-t^2})$ for $-1\le t\le 1$ and $=0$ otherwise. $f$ is obviously $C^\infty$ and compactly supported in the interval $[-1,1]$. We then construct the function

$$
\psi(u)=\frac{\int_{-1}^{u} f(t)\,dt}{\int_{-1}^{1} f(t)\,dt}.
$$

$\psi$ is $C^\infty$, nondecreasing and s.t. $\psi(-1)=0$, $\psi(1)=1$. Then the function $\varphi$ is obtained easily by joining 0 and 1 with $\psi$ suitably rescaled. Remark that in practice one needs only to compute the function $b$ at the points $\frac{\ell}{B^j}$. Therefore, once the maximal-resolution is known, these values can be computed and stored once for all. An instance of a needlet function is given in Figure 2.

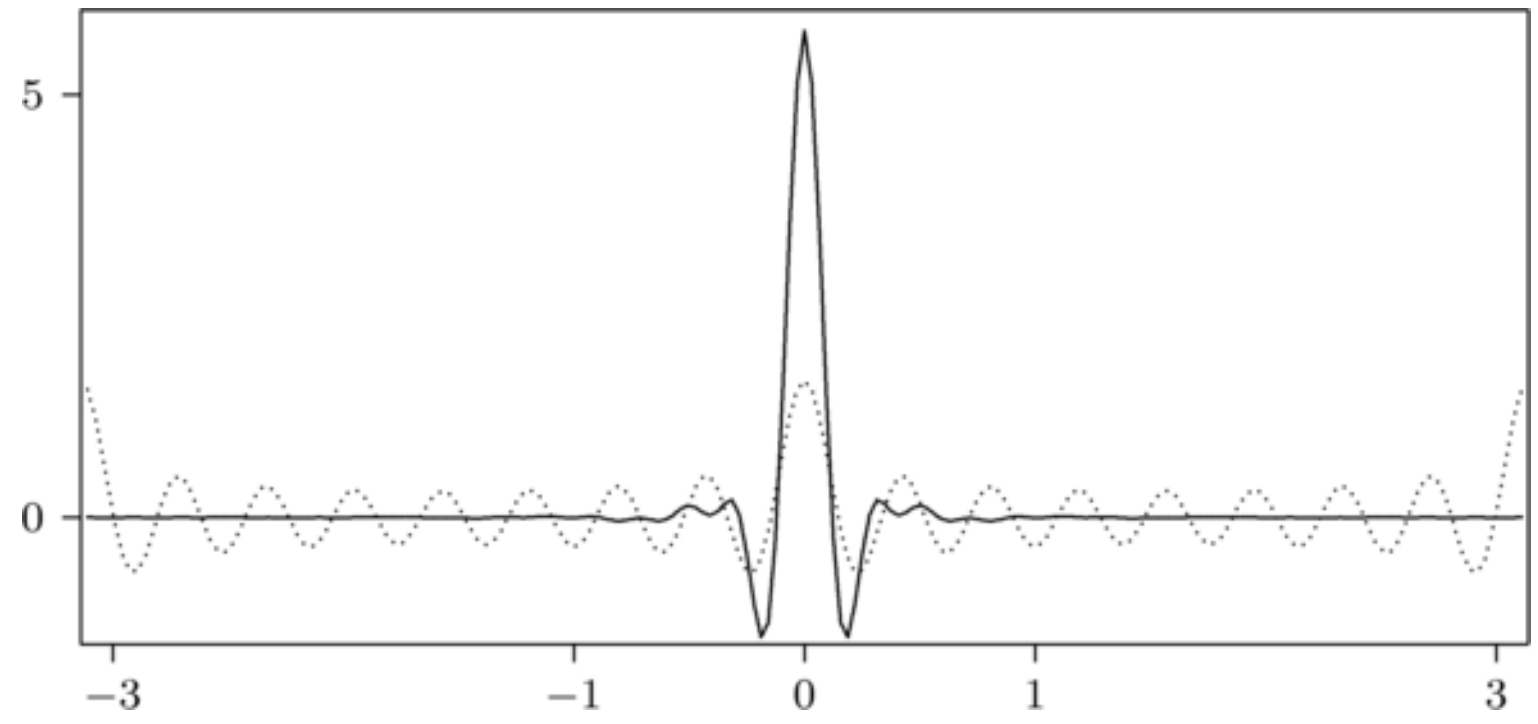

FIG. 2. *Typical graph of the needlet and a corresponding spherical harmonic (dots) as functions of the geodesic distance. Here $j=5$ and $B=2$. The localizing effect of the Littlewood–Paley device is remarkable.*

The random needlet coefficients are now evaluated as

$$\begin{aligned}\beta_{jk} &= \int_{\mathbb{S}^2} T(x)\psi_{jk}(x)\,dx \\ &= \sqrt{\lambda_{jk}}\sum_{\ell} b\left(\frac{\ell}{B^j}\right)\sum_{m=-\ell}^{\ell}\left\{\int_{\mathbb{S}^2} T(x)\overline{Y}_{\ell m}(x)\,dx\right\}Y_{\ell m}(\xi_{jk}) \qquad (16)\\ &= \sqrt{\lambda_{jk}}\sum_{\ell} b\left(\frac{\ell}{B^j}\right)\sum_{m=-\ell}^{\ell} a_{\ell m}Y_{\ell m}(\xi_{jk}).\end{aligned}$$

The practical implementation of (16) on a given random field requires the evaluation of its spherical harmonic coefficients $(a_{\ell m})$. In principle, the latter can be recovered by means of (9). In practice, in applications such as CMB data analysis the random field is continuously observed by means of antennae which average observations over tiny equal-area regions covering the whole sky; the resulting values are projected on a discretized grid, where the locations of points in the grid are chosen in order to make possible the approximation of (9) by means of cubature formulae; a standard package for this routine is HealPix, described in Górski et al. (2005). The final output of this algorithm is indeed a triangular array of coefficients $(a_{\ell m})$, but one may wonder whether numerical approximations may indeed spoil the validity of the theoretical results presented in the previous sections.

To investigate this claim, we produce some numerical evidence on one of the key properties of random needlet coefficients, that is, the uncorrelation across different scales (Figure 3). We simulated 100 independent copies of a random field, using the expansion (8). The coefficients $a_{lm}$ were sampled as independent [but for the condition $(-1)^m \overline{a}_{lm} = a_{l,-m}$] complex Gaussian r.v.'s with variance $C_l$. The results look encouraging: the actual correlation for all $j - j' \geq 2$ is in the order of 0.1–1%, which is indeed consistent with theoretical predictions, up to minor rounding errors.

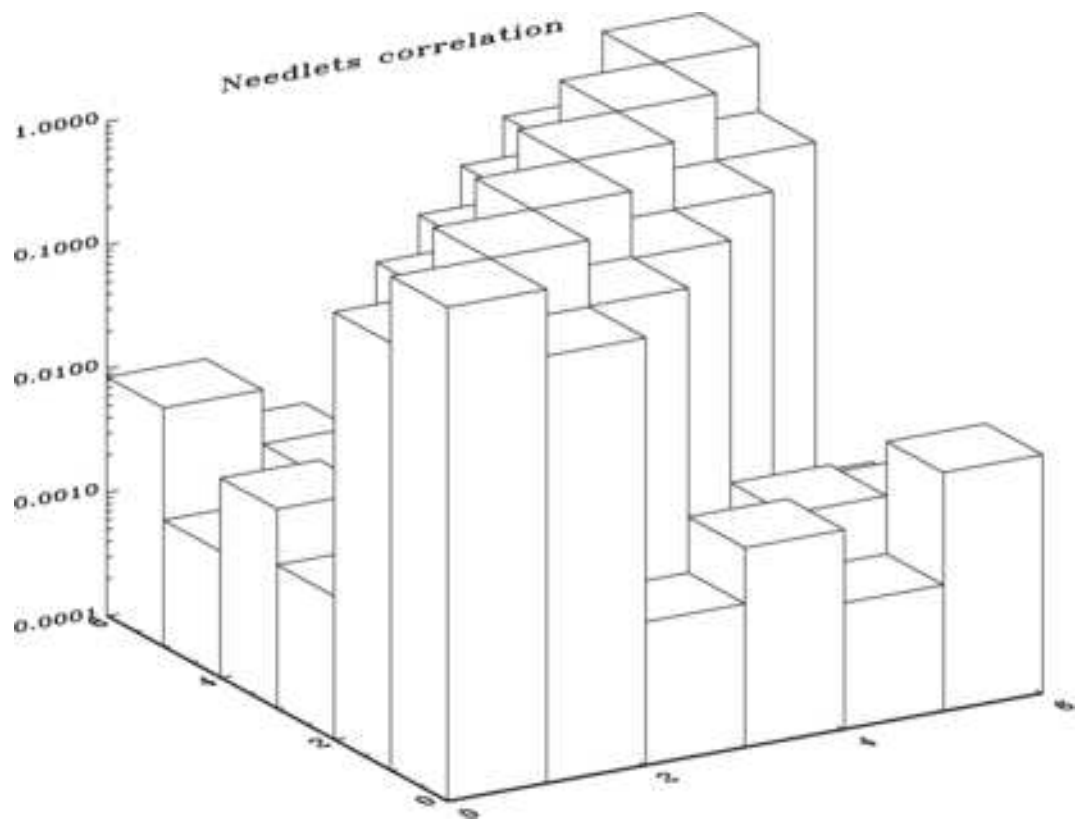

FIG. 3. *Decay of correlation on a CMB-like map. $j$ and $j'$ are on the two axes.*

We also performed some Monte Carlo experiments on the effect of missing observations on the values of the needlet coefficients. More precisely, for different types of sky gaps, we provide estimates of the quantity

$$D_{jk} := \frac{\mathrm{E}[\beta_{jk} - \widetilde{\beta}_{jk}]^2}{\mathrm{E}[\beta_{jk}^2]}. \tag{17}$$

First we mimicked the experimental data on the CMB radiation, as described for instance by the WMAP team (see http://map.gsfc.nasa.gov/). In particular, data on CMB are contaminated mainly by the presence of the Milky Way (which is located around the equator, in the standard choice of coordinates) and several so-called point sources, amounting basically to known clusters of galaxies which produce a radiation unrelated with CMB. To remove these emissions, the WMAP team has set to 0 the value of the field in a certain region, which is known as the Kp0 mask.

We simulated again 100 independent copies of a random field. The function $C_l$ was chosen in order to mimick the best fit from satellite observations of CMB [see Pietrobon, Balbi and Marinucci (2006) for details]. We fixed $B = 1.5$ and $j = 11$, corresponding to a range of frequencies from $l = 58$ to $l = 129$. We then estimated both the needlet coefficients $\widetilde{\beta}_{jk}$ (in the presence of missing observations) and $\beta_{jk}$ (for the completely observed field) and evaluated the gap between the two using the discrepancy $D_{jk}$ of (17).

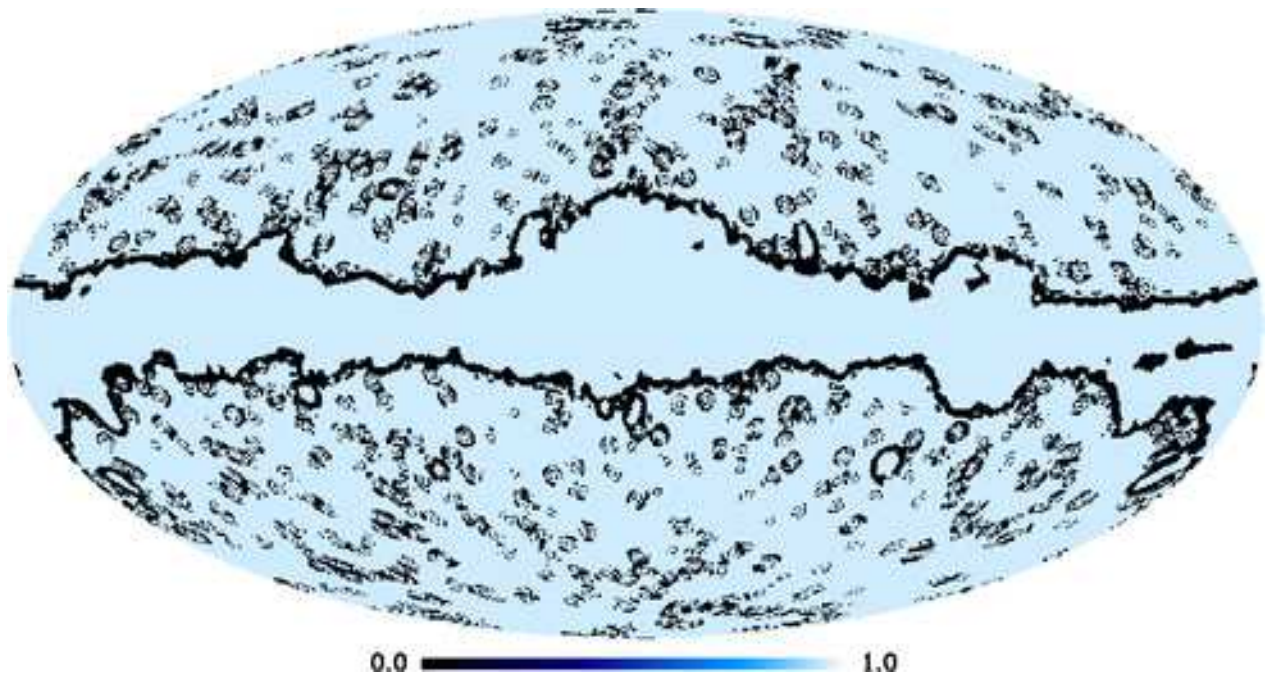

FIG. 4. *Localization properties on a CMB-like map.*

The results are displayed in Figure 4, where directions corresponding to a value of $D_{jk} > 0.1$ are marked with a black dot. Note that, even for such small values of $j$, the difference between $\beta_{jk}$ and $\widetilde{\beta}_{jk}$ is rather small; indeed $D_{jk}$ is above the threshold in approximately 20% of the cubature points. As expected, these points cluster in the neighborhoods of the mask. Refer to Guilloux, Fay and Cardoso (2007) for further numerical evidence.

**Acknowledgment.** We are grateful to D. Pietrobon for providing the numerical results in Section 8.

P. Baldi
D. Marinucci
Dipartimento di Matematica
Università di Roma Tor Vergata
Via della Ricerca Scientifica
00161 Roma
Italy
E-mail: baldi@mat.uniroma2.it
marinucc@mat.uniroma2.it

G. Kerkyacharian
D. Picard
Laboratoire de Probabilités
et Modèles Aléatoires
2, Pl. Jussieu
75251 Paris Cedex 05
France
E-mail: kerk@math.jussieu.fr
picard@math.jussieu.fr